\documentclass{article}

\usepackage{proof}
\usepackage{latexsym}
\usepackage{amssymb}

\usepackage[all]{xy}

\newcommand{\blem}{\begin{lemma}}
\newcommand{\elem}{\end{lemma}}
\newcommand{\bth}{\begin{theorem}}
\newcommand{\ethm}{\end{theorem}}
\newcommand{\benu}{\begin{enumerate}}
\newcommand{\eenu}{\end{enumerate}}
\newcommand{\bdes}{\begin{description}}
\newcommand{\edes}{\end{description}}
\newcommand{\bdf}{\begin{definition}}
\newcommand{\edf}{\end{definition}}
\newcommand{\bcor}{\begin{cor}}
\newcommand{\ecor}{\end{cor}}
\newcommand{\bprp}{\begin{proposition}}
\newcommand{\eprp}{\end{proposition}}
\newcommand{\bmlem}{\begin{mlemma}}
\newcommand{\emlem}{\end{mlemma}}
\newcommand{\bclm}{\begin{claim}}
\newcommand{\eclm}{\end{claim}}
\newcommand{\bprf}{{\bf Proof}.\hspace{2mm}}

\newcommand{\eprf}{\hspace*{\fill} $\Box$}

\newcommand{\beqn}{\begin{equation}}
\newcommand{\eeqn}{\end{equation}}
\newcommand{\beqnarr}{\begin{eqnarray}}
\newcommand{\eeqnarr}{\end{eqnarray}}
\newcommand{\beqnarrs}{\begin{eqnarray*}}
\newcommand{\eeqnarrs}{\end{eqnarray*}}

\newcommand{\spand}{\,\&\,}

\newcommand{\bfG}{\mbox{\boldmath$G$} }

\newcommand{\Natural}{\mathbb{N}}

\newcommand{\restrict}{\!\upharpoonright\!}

\newtheorem{theorem}{Theorem}[section]
\newtheorem{definition}[theorem]{Definition}
\newtheorem{proposition}[theorem]{Proposition}
\newtheorem{lemma}[theorem]{Lemma}
\newtheorem{cor}[theorem]{Corollary}
\newtheorem{mlemma}[theorem]{Main Lemma}
\newtheorem{claim}[theorem]{Claim}

\newcommand{\alp}{\alpha}

\newcommand{\veps}{\varepsilon}
\newcommand{\del}{\delta}
\newcommand{\Del}{\Delta}
\newcommand{\ome}{\omega}

\newcommand{\bet}{\beta}
\newcommand{\gam}{\gamma}
\newcommand{\Gam}{\Gamma}

\newcommand{\Sig}{\Sigma}

\newcommand{\lam}{\lambda}
\newcommand{\Lam}{\Lambda}
\newcommand{\vphi}{\varphi}

\newcommand{\fal}{\forall}
\newcommand{\exi}{\exists}

\newcommand{\Rarw }{\Rightarrow}

\newcommand{\lrarw}{\leftrightarrow}
\newcommand{\Lrarw}{\Leftrightarrow}

\newcommand{\calt}{{\cal T}}

\newcommand{\la}{\langle}
\newcommand{\ra}{\rangle}

\newcommand{\msfiv}{\mbox{\hspace{5mm}}}

\newcommand{\brem}{\begin{remark}}
\newcommand{\erem}{\end{remark}}
\newtheorem{remark}[theorem]{Remark}

\title{Derivatives of normal functions and $\ome$-models
}
\author{
Toshiyasu Arai
\\
Graduate School of Science,
Chiba University
\\
1-33, Yayoi-cho, Inage-ku,
Chiba, 263-8522, JAPAN
\\
tosarai@faculty.chiba-u.jp
}
\date{}

\begin{document}
\maketitle
\begin{abstract}
In this note the well-ordering principle for the derivative ${\sf g}^{\prime}$ 
of normal functions ${\sf g}$ on ordinals
is shown to be equivalent to the existence of arbitrarily large countable coded $\ome$-models of
the well-ordering principle for the function ${\sf g}$.
\end{abstract}

\section{Well-ordering principles}

In this note
 we are concerned with a proof-theoretic strength of a $\Pi^{1}_{2}$-statement ${\rm WOP}({\sf g})$ 
saying that
`for any well-ordering $X$, ${\sf g}(X)$ is a well-ordering',
where ${\sf g}: \mathcal{P}(\Natural)\to\mathcal{P}(\Natural)$ is a computable functional on sets $X$ of natural numbers.
$\la n,m\ra$ denotes an elementary recursive pairing function on $\Natural$.

\bdf
{\rm
$X\subset\Natural$ defines a binary relation 
$<_{X}:=\{(n,m): \la n,m\ra\in X\}$.
\beqnarrs
{\rm LO}(X) & :\Lrarw &
[
\fal n(n\not<_{X}n)\land \fal n,m,k(n<_{X}m<_{X}k\to n<_{X}k)
\\
&&
\land 
\fal n,m(n<_{X}m\lor n=m\lor m<_{X}n)
]
\\
{\rm Prg}[<_{X},Y] & :\Lrarw & \fal m\left(\fal n<_{X}m\,Y(n)\to Y(m)\right)
\\
{\rm TI}(<_{X},Y) & :\Lrarw & {\rm Prg}[<_{X},Y] \to\fal n\, Y(n)
\\
{\rm WO}(X) & :\Lrarw & {\rm LO}(X) \land \fal Y\, {\rm TI}(<_{X},Y)
\eeqnarrs

For a functional ${\sf g}: \mathcal{P}(\Natural)\to\mathcal{P}(\Natural)$,
\[
{\rm WOP}({\sf g}) :\Lrarw \fal X\left( {\rm WO}(X)\to {\rm WO}({\sf g}(X)) \right)
\]
}
\edf

The theorem due to J.-Y. Girard is a base for further results on the 
strengths of the well-ordering principles ${\rm WOP}({\sf g})$.

\bth\label{th:Girardp.299}{\rm (Girard\cite{Girard}, also cf.\,\cite{Hirst})}
\\
Over ${\rm RCA}_{0}$, ${\rm ACA}_{0}$ is equivalent to ${\rm WOP}(\lam X.\ome^{X})$.
\end{theorem}

The following theorem summarizes some known results on the strengths of ${\rm WOP}({\sf g})$ for 
${\sf g}$ larger than the exponential function.
${\rm ACA}_{0}^{+}$ is an extension of ${\rm ACA}_{0}$ by the axiom of the existence of
the $\ome$-th jump of a given set.
$\vphi\alp\bet=\vphi_{\alp}(\bet)$ denotes the binary Veblen function starting with $\ome^{\alp}$.

\bth\label{th:WOP}
\benu
\item\label{th:MontalbanMarcone}{\rm (Marcone and Montalb\'an\cite{MontalbanMarcone})}
\\
Over ${\rm RCA}_{0}$, ${\rm ACA}_{0}^{+}$ is equivalent to ${\rm WOP}(\lam X.\veps_{X})$.
\item\label{th:Friedman}{\rm (H. Friedman)}
\\
Over ${\rm RCA}_{0}$, ${\rm ATR}_{0}$ is equivalent to ${\rm WOP}(\lam X.\vphi X0)$.
\eenu
\end{theorem}

Theorem \ref{th:WOP} is proved
in \cite{MontalbanMarcone} computability theoretically.
M. Rathjen noticed that the principle ${\rm WOP}({\sf g})$ is tied to the existence of
\textit{countable coded $\ome$-models}.

\bdf
{\rm
A \textit{countable coed $\ome$-model} of a second-order arithmetic $T$ is a set $Q\subset\Natural$ such that
$M(Q)\models T$, where
$M(Q)=\la\Natural,\{(Q)_{n}\}_{n\in\Natural}, +,\cdot,0,1,<\ra$
with $(Q)_{n}=\{m\in\Natural: \la n,m\ra\in Q\}$.

Let $X\in_{\ome}Y:\Lrarw(\exi n[X=(Y)_{n}])$ and
$X=_{\ome}Y:\Lrarw(\fal Z(Z\in_{\ome}X\lrarw Z\in_{\ome}Y))$.
}
\edf
It is not hard to see that over ${\rm ACA}_{0}$,
${\rm ACA}_{0}^{+}$ is equivalent to
the fact that there exists an arbitrarily large countable coded $\ome$-model of ${\rm ACA}_{0}$, cf.\,\cite{AfshariRathjen} and Lemma \ref{lem:Veblenjumphier} below.
The fact means that there is a countable coded $\ome$-model $Q$ of ${\rm ACA}_{0}$
containing a given set $X$, i.e., $X=(Q)_{0}$.
From this characterization, Afshari and Rathjen\cite{AfshariRathjen} gives 
a purely proof-theoretic proof of Theorem \ref{th:WOP}.\ref{th:MontalbanMarcone}.
Their proof is based on Sch\"utte's method of complete proof search in $\ome$-logic.
The proof is extended by Rathjen and Weiermann\cite{RathjenWeiermannRM}
to give an alternative proof of Theorem \ref{th:WOP}.\ref{th:Friedman}.
Furthermore Rathjen\cite{RathjenHFriedman} lifts Theorem \ref{th:WOP}.\ref{th:MontalbanMarcone} up to $\Gam$-function and ${\rm ATR}_{0}$ as follows.

\bdf
{\rm
A continuous and strictly increasing function on ordinals is said to be a \textit{normal} function.

For a normal function $f$, its \textit{derivative} $f^{\prime}$ is a normal function
enumerating the fixed points of the function $f$.
}
\edf
The $(\alp+1)$-th branch $\vphi_{\alp+1}:\bet\mapsto\vphi_{\alp+1}(\bet)$ of the Veblen function
is the derivative $(\vphi_{\alp})^{\prime}$ of the previous one $\vphi_{\alp}$, and for limit $\lam$, 
$\vphi_{\lam}$ enumerates
the common fixed points of the functions $\vphi_{\alp}\,(\alp<\lam)$.
The $\Gam$-function $\alp\mapsto\Gam_{\alp}$ is the derivative of the normal function
$\alp\mapsto\vphi\alp 0$.

\bth\label{th:RathjenHFriedman}
{\rm (Rathjen\cite{RathjenHFriedman})}
\\
Over ${\rm RCA}_{0}$, ${\rm WOP}(\lam X.\Gam_{X})$ is equivalent to
the existence of arbitrarily large countable coded $\ome$-models of ${\rm ATR}_{0}$.
\end{theorem}

In view of Theorem \ref{th:Girardp.299},
Theorem \ref{th:WOP}.\ref{th:MontalbanMarcone}
is equivalently stated: over ${\rm RCA}_{0}$, 
${\rm WOP}(\lam X.\veps_{X})$ is equivalent to 
the existence of arbitrarily large countable coded $\ome$-models of ${\rm WOP}(\lam X.\ome^{X})$.
Moreover relying on \ref{th:WOP}.\ref{th:Friedman},
Theorem \ref{th:RathjenHFriedman} states that over ${\rm RCA}_{0}$, 
${\rm WOP}(\lam X.\Gam_{X})$ is equivalent to
the existence of arbitrarily large countable coded $\ome$-models of ${\rm WOP}(\lam X.\vphi X 0)$.
Here is a striking similarity: 
$\lam\alp.\veps_{\alp}$ is the derivative of the function $\lam\alp.\ome^{\alp}$,
and $\lam\alp.\Gam_{\alp}$ is the one of $\lam\alp.\vphi\alp 0$.

\bdf
{\rm $T^{+}$ denotes the extension of a second-order arithmetic $T$ by the axiom stating that
\beqn\label{eq:aboveT}
\mbox{
`there exists an arbitrarily large countable coded $\ome$-model of $T$'}
\eeqn
}
\edf
Note that when $T$ is axiomatized by a $\Pi^{1}_{2}$-sentence over ${\rm RCA}_{0}$,
$T^{+}$ is axiomatized by the $\Pi^{1}_{2}$-sentence (\ref{eq:aboveT}) over ${\rm RCA}_{0}$.

These results suggest us a general fact:
\beqn\label{quote:1}
\mbox{${\rm WOP}({\sf g}^{\prime})$ is equivalent to ${\rm WOP}({\sf g})^{+}$ over ${\rm ACA}_{0}$.}
\eeqn
In this note we confirm it for a variety of normal functions ${\sf g}$.
Theorem \ref{th:WOP}.\ref{th:MontalbanMarcone} follows 
from (\ref{quote:1}) for ${\sf g}(\alp)=\ome^{\alp}$, and
Theorem \ref{th:RathjenHFriedman} from 
Theorem \ref{th:WOP}.\ref{th:Friedman} and (\ref{quote:1}) for ${\sf g}(\alp)=\vphi_{\alp}(0)$.

We assume that the normal function ${\sf g}$ enjoys the following conditions.
The computability of the functional ${\sf g}$ and the linearity of ${\sf g}(X)$ for linear orderings $X$ are assumed.
Moreover
${\sf g}(X)$ is assumed to be a \textit{term structure} over
constants ${\sf g}(c)\, (c\in X)$ and some function symbols $f$.
For the term structures ${\sf G}(X)=\left({\sf g}(X),<_{{\sf g}(X)};f,\ldots\right)$ we need two facts:
First 
 if $\left( X,<_{X}\right)$ is a substructure of $\left( Y,<_{Y}\right)$, then 
${\sf G}(X)$ is a substructure of ${\sf G}(Y)$.
Second $\la {\sf g}(c): c\in X\ra$ is an indiscernible sequence for ${\sf G}(X)$.
These two postulates allow us to extend 
an order preserving map $f$ between linear orderings $X,Y$ to 
an order preserving map $F$ between ${\sf g}(X)$ and ${\sf g}(Y)$, cf.\,Proposition \ref{prp:dilate}:
\[
\xymatrix{
{\sf g}(X) \ar[r]^{F}  & {\sf g}(Y) 
\\
X \ar[r]_{f} \ar[u]^{i} & Y \ar[u]_{i}
}
\]
Moreover we assume that
$\left({\sf g}^{\prime}(X); 0,+,\lam \alp.\ome^{\alp}\right)$ is a substructure of 
the term structure ${\sf G}^{\prime}(X)$ for the derivative ${\sf g}^{\prime}$.
Then (\ref{quote:1}) is shown in Theorem \ref{th:derivativemodel}.

Next (\ref{quote:1}) suggests us a result on common fixed points.
Let $\vphi[{\sf g}]_{\alp}(\bet)$ denote the $\alp$-th Veblen function starting with 
$\vphi[{\sf g}]_{0}(\bet)={\sf g}(\bet)$.
For $\alp>0$
\beqn\label{quote:2}
\mbox{${\rm WOP}(\vphi[{\sf g}]_{\alp})$ is equivalent to
$\left(\fal\bet<\alp{\rm WOP}(\vphi[{\sf g}]_{\bet})\right)^{+}$ over ${\rm ACA}_{0}+{\rm LO}(\alp)$.}
\eeqn
Under a mild condition on the Veblen hierarchy $\{\vphi[{\sf g}]_{\alp}\}_{\alp}$,
we confirm (\ref{quote:2}) in Theorem \ref{th:commonderivative}.
\\

Next consider ${\rm WOP}(\vphi_{\alp})$ with 
$\vphi_{\alp}=\vphi[{\sf g}]_{\alp}$ for the most familiar ${\sf g}(\bet)=\ome^{\bet}$.

Let ${\rm TJ}(X)$ denote the Turing jump of sets $X$.
${\rm Hier}_{\alp}(X,Y)$ designates that $\{(Y)_{\bet}\}_{\bet<\alp}$ is the Turing jump hierarchy
starting with $X=(Y)_{0}$ for the least element $0$ in the ordering $<$: for any non-zero $\bet<\alp$,
if $\bet=\gam+1$, then $(Y)_{\bet}={\rm TJ}((Y)_{\gam})$, and when $\bet$ is limit,
$(Y)_{\bet}=\sum_{\gam<\bet}(Y)_{\gam}$.

The following Lemma \ref{lem:Veblenjumphier} is shown in \cite{MontalbanMarcone}, Theorem 1.9,
and it yields
Theorem \ref{th:WOP}.\ref{th:Friedman}.

\blem\label{lem:Veblenjumphier}
{\rm (\cite{MontalbanMarcone})}.

Over ${\rm ACA}_{0}+{\rm WO}(\alp)$,
${\rm WOP}(\vphi_{\alp})$ is equivalent to
$\fal X\exi Y\, {\rm Hier}_{\ome^{\alp}}(X,Y)$.
\elem
\bprf
It is well known that ${\rm WOP}(\vphi_{\alp})$ follows from
$\fal X\exi Y\, {\rm Hier}_{\ome^{\alp}}(X,Y)+{\rm WO}(\alp)$.

Let $A(\alp):\Lrarw[{\rm WOP}(\vphi_{\alp}) \to \fal X\exi Y\, {\rm Hier}_{\ome^{\alp}}(X,Y)]$.
It suffices to show in ${\rm ACA}_{0}$ that $A(\alp)$ assuming $A(\bet)$ holds for any $\bet<\alp$ 
in any countable coded $\ome$-models of ${\rm ACA}_{0}$.
Then ${\rm WO}(\alp)$ yields $A(\alp)$.

Assume that $A(\bet)$ holds for any $\bet<\alp$ 
in any countable coded $\ome$-models of ${\rm ACA}_{0}$.
Suppose ${\rm WOP}(\vphi_{\alp})$ for $\alp>0$.
Then by (\ref{quote:2}) we have $(\fal\bet<\alp{\rm WOP}(\vphi_{\bet}))^{+}$.
Given a set $X$, pick a countable coded $\ome$-model $Z$ of $\fal\bet<\alp{\rm WOP}(\vphi_{\bet})$
such that $X\in_{\ome}Z$.
$Z$ is an $\ome$-model of ${\rm ACA}_{0}$ by Theorem \ref{th:Girardp.299}).
By the assumption we obtain $\fal\bet<\alp\fal X\exi ! Y\, {\rm Hier}_{\ome^{\bet}}(X,Y)$ in $Z$.
Given a set $X$ let
$W=\{\la \gam,m\ra: \gam<\ome^{\alp}, Z\models\exi Y[{\rm Hier}_{\gam}(X,Y)\land m\in(Y)_{\gam}]\}$.
$W$ is a set by ${\rm ACA}_{0}$.
If $\alp$ is a limit number, then 
$\fal \bet<\alp {\rm Hier}_{\ome^{\bet}}(X,W)$ yields ${\rm Hier}_{\ome^{\alp}}(x,W)$.
When $\alp=\bet+1$, we see by induction on $k<\ome$ that
$\fal k<\ome {\rm Hier}_{\ome^{\bet}k}(X,W)$, and hence ${\rm Hier}_{\ome^{\alp}}(x,W)$.
\eprf

\section{Term structures}\label{sect:termstr}

Let us compare the proof-theoretic strength ${\rm WOP}({\sf g}^{\prime})$ with ${\rm WOP}({\sf g})$
for normal function ${\sf g}$.
First of all, both ${\sf g}^{\prime}$ and ${\sf g}$ need to be definable to express
formulas ${\rm WOP}({\sf g}^{\prime})$ and ${\rm WOP}({\sf g})$ in $\Pi^{1}_{2}$-formulas.
Moreover the fact that ${\sf g}$ sends
linear orderings $X$ to linear orderings ${\sf g}(X)$ should be provable in an elementary way.
However we need stronger conditions.

${\sf g}$ sends a binary relation $<_{X}$ on a set $X$ to a binary relation
$<_{{\sf g}(X)}={\sf g}(<_{X})$ on a set ${\sf g}(X)$.
We assume that ${\sf g}(X)$ is a Skolem hull, i.e., a term structure over
constants ${\sf g}(c)\, (c\in \{0\}\cup X)$ with the least element $0$ in the order $<_{X}$,
and some (possibly infinite number of) 
function symbols $f\in\mathcal{F}$.
Let us assume that each function symbol except $+$ has a fixed arity.
Some function symbol $f\in\mathcal{F}$ may not be totally defined.
In other words $f(\bet_{1},\ldots,\bet_{n})$ may be an illegal expression for $\bet_{1},\ldots,\bet_{n}\in{\sf g}(X)$,
i.e., $f(\bet_{1},\ldots,\bet_{n})\not\in{\sf g}(X)$.

\bdf\label{df:gtmstr}
{\rm
\benu
\item
${\sf g}(X)$ is said to be a \textit{computably linear} term structure if
there are three $\Sig^{0}_{1}(X)$-formulas ${\sf g}(X), <_{{\sf g}(X)},=$ 
for which all of the following facts are provable in ${\rm RCA}_{0}$:
let $\alp,\bet,\gam,\ldots$ range over terms.
\benu
\item(Computability)
Each of ${\sf g}(X)$, $<_{{\sf g}(X)}$ and $=$ is $\Del^{0}_{1}(X)$-definable.
${\sf g}(X)$ is a computable set, and $<_{{\sf g}(X)}$ and $=$ are computable binary relations.

\item(Congruence)
$=$ is a congruence relation on the structure $\la{\sf g}(X);<_{{\sf g}(X)},f,\ldots\ra$.

Let us denote ${\sf g}(X)/=$ the quotient set.

In what follows assume that $<_{X}$ is a linear ordering on $X$.

 \item(Linearity)
 $<_{{\sf g}(X)}$ is a linear ordering on ${\sf g}(X)/=$.






\item(Increasing)
${\sf g}$ is strictly increasing:
$c<_{X}d \Rarw {\sf g}(c)<_{{\sf g}(X)}{\sf g}(d)$.

\item(Continuity)
${\sf g}$ is continuous:
Let $\alp<_{{\sf g}(X)}{\sf g}(c)$ for a limit $c\in X$ and $\alp\in{\sf g}(X)$.
Then there exists a $d<_{X}c$ such that $\alp<_{{\sf g}(X)}{\sf g}(d)$.
\eenu

\item
A computably linear term structure ${\sf g}(X)$ is said to be \textit{extedible} if it enjoys
the following two conditions.
\benu
\item(Suborder)
If $\la X,<_{X}\ra$ is a substructure of $\la Y,<_{Y}\ra$, then 
$\la{\sf g}(X);=, <_{{\sf g}(X)}, f,\ldots\ra$ is a substructure of $\la{\sf g}(Y); =, <_{{\sf g}(Y)},f,\ldots\ra$.

\item(Indiscernible)
$\la {\sf g}(c): c\in\{0\}\cup X\ra$ is an indiscernible sequence for 
linera orderings $\la{\sf g}(X),<_{{\sf g}(X)}\ra$:
Let $\alp[{\sf g}(c_{1}),\ldots,{\sf g}(c_{n})], \bet[{\sf g}(c_{1}),\ldots,{\sf g}(c_{n})]\in{\sf g}(X)$ be
terms such that constants occurring in them are among the list ${\sf g}(c_{1}),\ldots,{\sf g}(c_{n})$.
Then for any increasing sequences $c_{1}<_{X}\ldots<_{X}c_{n}$ and
$d_{1}<_{X}\ldots<_{X}d_{n}$
\beqnarr
&&
\alp[{\sf g}(c_{1}),\ldots,{\sf g}(c_{n})]<_{{\sf g}(X)}\bet[{\sf g}(c_{1}),\ldots,{\sf g}(c_{n})]
\label{eq:indiscerng}
\\
& \Lrarw &
\alp[{\sf g}(d_{1}),\ldots,{\sf g}(d_{n})]<_{{\sf g}(X)}\bet[{\sf g}(d_{1}),\ldots,{\sf g}(d_{n})]
\nonumber
\eeqnarr

\eenu
\eenu
}
\edf

\bprp\label{prp:dilate}
{\rm Suppose ${\sf g}(X)$ is an extendible term structure.
Then the following is provable in ${\rm RCA}_{0}$:}
Let both $X$ and $Y$ be linear orderings.

Let $f:\{0\}\cup X\to \{0\}\cup Y$ be an order preserving map, $n<_{X}m\Rarw f(n)<_{Y}f(m)\, (n,m\in\{0\}\cup X)$.
Then there is an order preserving map $F:{\sf g}(X)\to{\sf g}(Y)$, $n<_{{\sf g}(X)}m\Rarw F(n)<_{{\sf g}(Y)}F(m)$.
\eprp
\bprf
Let $\alp[{\sf g}(c_{1}),\ldots,{\sf g}(c_{n})]\in{\sf g}(X)$ be
a term such that constants occurring in it are among the list ${\sf g}(c_{1}),\ldots,{\sf g}(c_{n})$ for 
$c_{i}\in\{0\}\cup X$.

Define
$F(\alp[{\sf g}(c_{1}),\ldots,{\sf g}(c_{n})])=\alp[{\sf g}(f(c_{1})),\ldots,{\sf g}(f(c_{n}))]$.
From (\ref{eq:indiscerng}) on ${\sf g}(X+Y)$,
we see that $F$ is an order preserving map from ${\sf g}(X)$ to ${\sf g}(Y)$.
Moreover we see that 
$\alp[{\sf g}(c_{1}),\ldots,{\sf g}(c_{n})]=\bet[{\sf g}(c_{1}),\ldots,{\sf g}(c_{n})] \Rarw
F(\alp[{\sf g}(c_{1}),\ldots,{\sf g}(c_{n})])=F(\bet[{\sf g}(c_{1}),\ldots,{\sf g}(c_{n})])$.
\eprf

\bdf\label{df:omegag}
{\rm Suppose that function symbols $+,\ome$ are in the list $\mathcal{F}$ of function symbols for
a computably linear term structure ${\sf g}(X)$. Let $1:=\ome^{0}$, and $2:=1+1$, etc.

${\sf g}(X)$ is said to be an \textit{exponential} term structure 
(with respect to function symbols $+,\ome$) if all of the followings are provable in ${\rm RCA}_{0}$.
\benu
 \item\label{df:omegag.1}
 $0$ is the least element in $<_{{\sf g}(X)}$, and $\alp+1$ is the successor of $\alp$.
 \item\label{df:omegag.2}
 $+$ and $\ome$ enjoy some familiar conditions.
 \benu
 \item
  $\alp<_{{\sf g}(X)}\bet\to \ome^{\alp}+\ome^{\bet}=\ome^{\bet}$.
  \item
 $\gam+\lam=\sup\{\gam+\bet:\bet<\lam\}$ when $\lam$ is a limit number, i.e.,
 $\lam\neq 0$ and $\fal \bet<_{{\sf g}(X)}\lam(\bet+1<_{{\sf g}(X)}\lam)$.
 \item
  $\bet_{1}<_{{\sf g}(X)}\bet_{2} \to \alp+\bet_{1}<_{{\sf g}(X)}\alp+\bet_{2}$, and
  $\alp_{1}<_{{\sf g}(X)}\alp_{2} \to \alp_{1}+\bet\leq_{{\sf g}(X)}\alp_{2}+\bet$.
  \item
  $(\alp+\bet)+\gam=\alp+(\bet+\gam)$.
  \item
  $\alp<_{{\sf g}(X)}\bet\to\exi\gam\leq_{{\sf g}(X)}\bet(\alp+\gam=\bet)$.
  \item
  Let $\alp_{n}\leq_{{\sf g}(X)}\cdots\leq_{{\sf g}(X)}\alp_{0}$ and 
  $\bet_{m}\leq_{{\sf g}(X)}\cdots\leq_{{\sf g}(X)}\bet_{0}$.
  Then $\ome^{\alp_{0}}+\cdots+\ome^{\alp_{n}}<_{{\sf g}(X)}\ome^{\bet_{0}}+\cdots+\ome^{\bet_{m}}$ iff
  either $n<m$ and $\fal i\leq n(\alp_{i}=\bet_{i})$, or 
  $\exi j\leq\min\{n,m\}[\alp_{j}<_{{\sf g}(X)}\bet_{j}\land \fal i<j(\alp_{i}=\bet_{i})]$.
 \eenu

\item\label{df:omegag.3}
Each $f(\bet_{1},\ldots,\bet_{n})\in{\sf g}(X)\, (f\in\mathcal{F})$ as well as 
${\sf g}(c)\,(c\in\{0\}\cup X)$ 
is closed under $+$.
In other words the terms $f(\bet_{1},\ldots,\bet_{n})$ and ${\sf g}(c)$ denote additively closed ordinals
 (additive principal numbers)
when $<_{{\sf g}(X)}$ is a well ordering.

\eenu
}
\edf

In what follows we assume that ${\sf g}(X)$ is an extendible term structure,
and ${\sf g}^{\prime}(X)$ is an exponential term structure.
Constants in the term structure ${\sf g}^{\prime}(X)$ are ${\sf g}^{\prime}(c)$ for $c\in\{0\}\cup X$,
and function symbols in $\mathcal{F}\cup\{0,+\}\cup\{{\sf g}\}$ with a unary function symbol ${\sf g}$.
When $\mathcal{F}=\emptyset$, let $\ome^{\alp}:={\sf g}(\alp)$.
Otherwise we assume that $\ome$ is in the list $\mathcal{F}$.
Furthermore assume that ${\rm RCA}_{0}$ proves that
\beqnarr
\bet_{1},\ldots,\bet_{n}<_{{\sf g}^{\prime}(X)}{\sf g}^{\prime}(c) & \to &
 f(\bet_{1},\ldots,\bet_{n})<_{{\sf g}^{\prime}(X)}{\sf g}^{\prime}(c) \,(f\in\mathcal{F}\cup\{+,{\sf g}\})
 \nonumber
 \\
\ome^{{\sf g}^{\prime}(\bet)} & = & {\sf g}({\sf g}^{\prime}(\bet))={\sf g}^{\prime}(\bet)
\nonumber
\\
{\sf g}^{\prime}(0) & = & \sup_{n}{\sf g}^{n}(0)
\label{df:omegag.4}
\\
{\sf g}^{\prime}(c+1) & = & \sup_{n}{\sf g}^{n}({\sf g}^{\prime}(c)+1)\, (c\in\{0\}\cup X)
\nonumber
\eeqnarr
where ${\sf g}^{n}$ denotes the $n$-th iterate of the function ${\sf g}$,
and we are assuming in the last  that
the successor element $c+1$ of $c$ in $X$ exists.
Note that the last two in (\ref{df:omegag.4}) are true for normal functions ${\sf g}$ when ${\sf g}(0)>0$.

Assume that $<_{X}$ is a linear ordering.
Each non-zero term $\bet\in{\sf g}^{\prime}(X)$ is written as a Cantor normal form
$\bet=\bet_{1}+\cdots+\bet_{n}$ where $\bet_{n}\leq_{{\sf g}^{\prime}(X)}\ldots\leq_{{\sf g}^{\prime}(X)}\bet_{1}$ and
each $\bet_{i}$ is an $f$-term $f(\gam_{1},\ldots,\gam_{m})$ with $f\in\mathcal{F}$ or
${\sf g}^{\prime}(c)$.
Using the Cantor normal form, we can define the natural (commutative) sum $\alp\#\bet$ of
terms $\alp,\bet\in{\sf g}^{\prime}(X)$ which enjoys $\alp\#\bet=\bet\#\alp$ and
$\alp_{1}<_{{\sf g}^{\prime}(X)}\alp_{2} \Rarw \alp_{1}\#\bet<_{{\sf g}^{\prime}(X)}\alp_{2}\#\bet$.

\bth\label{th:derivativemodel}
Let ${\sf g}(X)$ be an extendible term structure, and 
${\sf g}^{\prime}(X)$ an exponential term structure for which (\ref{df:omegag.4}) holds.

Then the following two are mutually equivalent over ${\rm ACA}_{0}$:
\benu
\item
${\rm WOP}({\sf g}^{\prime})$.
\item
$\left({\rm WOP}({\sf g}^{\prime})\right)^{+}
:\Lrarw \fal X\exi Y[X\in Y \land M_{Y}\models{\rm WOP}({\sf g})]$.
Namely there exists an arbitrarily large countable coded $\ome$-model of ${\rm WOP}({\sf g})$.
\eenu
\end{theorem}

First let us show the easy half.
Let sets $X,U$ be given such that ${\rm WO}(<_{0})$ for $<_{0}=<_{X}$.
We have ${\rm LO}(<_{{\sf g}^{\prime}(X)})$.
Pick a countable coded $\ome$-model $M$ of ${\rm WOP}({\sf g})$ such that $X,U\in M$.
Then ${\sf g}(X),{\sf g}^{\prime}(X)\in M$.
Let $<_{1}$ be obtained from $<_{0}$ by adding the largest element $\alp$.
This means that $a<_{1}\alp$ for any $a$ in the field of $<_{0}$.
We have ${\rm WO}(<_{1})$ by ${\rm WO}(<_{0})$.
We show $Prg[<_{1}, C(a)]$ for an arithmetical formula
\[
C(a):\Lrarw M\models\fal Y\left( Prg[<_{2},Y]\to\fal x<_{2}{\sf g}^{\prime}(a)\, Y(x)\right)
\]
for $<_{2}=<_{{\sf g}^{\prime}(<_{1})}$.
This yields $C(\alp)$.
Since by (\ref{df:omegag.4}),
 $x<_{2}{\sf g}^{\prime}(\alp)$ for any $x$ in the field of $<_{{\sf g}^{\prime}(X)}=<_{{\sf g}^{\prime}(<_{0})}$,
we obtain $M\models\fal Y\left( Prg[<_{{\sf g}^{\prime}(X)},Y]\to\fal x\in {\rm fld}({\sf g}^{\prime}(X))\, Y(x)\right)$.
Hence we obtain
$M\models\left( Prg[<_{{\sf g}^{\prime}(X)},U]\to\fal x\in {\rm fld}({\sf g}^{\prime}(X))\, U(x)\right)$, i.e.,
${\rm TI}(<_{{\sf g}^{\prime}(X)},U)$.
Since $U$ is an arbitrary set, we conclude ${\rm WO}(<_{{\sf g}^{\prime}(X)})$.

It remains to show that $Prg[<_{1}, C(a)]$.
When $a$ is a limit element, this follows from the continuity of the function ${\sf g}^{\prime}(a)$.
Assuming $C(a)$, let us show $C(a+1)$.
Argue in the model $M$.
Suppose $Prg[<_{2},Y]$ and $x<_{2}{\sf g}^{\prime}(a+1)=\sup_{n}{\sf g}^{n}({\sf g}^{\prime}(a)+1)$ 
by (\ref{df:omegag.4}).
By induction on $n<\ome$ we see that $\fal x<_{2}{\sf g}^{n}({\sf g}^{\prime}(a)+1)\, Y(x)$
using ${\rm WOP}({\sf g})$ and $C(a)$, i.e., ${\rm WO}(<_{2}\restrict({\sf g}^{\prime}(a)+1))$.
Hence we obtain $C(a+1)$.
$C(0)$ is seen similarly.

\section{Proof search}

Conversely assume ${\rm WOP}({\sf g}^{\prime})$.
We need to find a countable coded $\ome$-model of ${\rm WOP}({\sf g})$.
The idea in \cite{AfshariRathjen,RathjenHFriedman,RathjenWeiermannRM} is
 to search a derivation of the negation of ${\rm WOP}({\sf g})$ in $\ome$-logic.
Construct a locally correct $\ome$-branching tree in a canonical way.
If the search results in a fail, i.e., if the constructed tree is not well-founded, 
 then we can believe in the consistency of ${\rm WOP}({\sf g})$ in $\ome$-logic.
In fact we can find a countable coded $\ome$-model of ${\rm WOP}({\sf g})$ from an infinite
path through the tree.
Otherwise the tree is well-founded, i.e., a derivation in a depth $\alp$.
It turns out that the derivation can be converted to a cut-free deduction with the empty sequent at its root,
and in depth bounded by ${\sf g}^{\prime}(\alp)$.
Then by our assumption ${\rm WOP}({\sf g}^{\prime})$, the deduction is well-founded, i.e., 
a derivation of the empty sequent.
We see that this is not the case by transfinite induction up to ${\sf g}^{\prime}(\alp)$.
This shows the consistency of ${\rm WOP}({\sf g})$ in $\ome$-logic based on 
${\rm WOP}({\sf g}^{\prime})$.
Now details follows.
\\

Let $\mathcal{Q}\subset\Natural$ be a given set, which is viewed as a family $\{(\mathcal{Q})_{i}: i<\ome\}$ 
of sets of natural numbers.
The language $\mathcal{L}_{\ome}$ here consists of function symbols for elementary recursive functions including
$0$ and the successor $S$, predicate symbols $=,\neq$ and
unary predicate variables $\{X_{i},E_{i}: i<\ome\}$ and their compliments $\bar{X}_{i},\bar{E}_{i}$.
Let us write $n<_{i}m$ for $n<_{X_{i}}m$, i.e., for $X_{i}(\la n,m\ra)$, and
$n<_{{\sf g}_{i}}m$ for $n<_{{\sf g}(X_{i})}m$.
Each $E_{i}$ is a fresh variable expressing the well foundedness ${\rm TI}(<_{i},E_{i})$ of the relation $<_{i}$.
Recall that each closed term $t$ is identified with its value $t^{\Natural}$, a numeral.

\[
{\rm D}_{\mathcal{Q}}(i,n)=\left\{
\begin{array}{ll}
X_{i}(n) & n\in(\mathcal{Q})_{i}
\\
\bar{X}_{i}(n) & n\not\in(\mathcal{Q})_{i}
\end{array}
\right.
\]
and ${\rm Diag}(\mathcal{Q})=\{{\rm D}_{\mathcal{Q}}(i,n): i,n\in\Natural\}$.

A \textit{true literal} is one of the form 
$t_{0}=t_{1}\, (t_{0}^{\Natural}=t_{1}^{\Natural})$,
$s_{0}\neq s_{1}\, (s_{0}^{\Natural}\neq s_{1}^{\Natural})$,
and ${\rm D}_{\mathcal{Q}}(i,n)$ for $i,n<\ome$.
\\
{\bf Axioms} in $\bfG(\mathcal{Q})+(prg)+(W)$ are
\[
\Gam,\bar{E}_{i}(n),E_{i}(n)
\]
and for true literals $L$
\[
\Gam,L
\]
{\bf Inference rules} are in $\bfG(\mathcal{Q})+(prg)+(W)$.
\[
\infer[(\lor)]{\Gam,A_{0}\lor A_{1}}
{\Gam,A_{0}\lor A_{1},A_{i}}
\,
\infer[(\land)]{\Gam,A_{0}\land A_{1}}
{
\Gam,A_{0}
&
\Gam,A_{1}
}
\,
\infer[(\exi)]{\Gam,\exi x\, A(x)}
{\Gam,\exi x\, A(x),A(n)}
\,
\infer[(\fal\ome)]{\Gam,\fal x\, A(x)}
{
\{\Gam,A(n): n<\ome\}
}
\]
$(\exi^{2})$ for $i<\ome$ and $(\fal^{2})$ with an eigenvariable $Z$
\[
\infer[(\exi^{2})]{\Gam,\exi Y\, A(Y)}
{\Gam,\exi Y\, A(Y), A(X_{i})
}
\msfiv
\infer[(\fal^{2})]{\Gam,\fal Y\, A(Y)}
{\Gam,A(Z)
}
\]
and the following two for $i,m<\ome$:
\[
\infer[(prg)_{i}]{\Gam, E_{i}(m)}
{
\{\Gam,E_{i}(n): n<_{i}m \mbox{ is true}\}
}
\]
where by saying that $n<_{i}m$ is true we mean $\la n,m\ra\in(\mathcal{Q})_{i}$.
\[
\infer[(W)_{i}]{\Gam}
{
\Gam,{\rm LO}(<_{i})
&
\Gam,\fal x\, E_{i}(x)
&
\exi Y\lnot {\rm TI}(<_{{\sf g}_{i}},Y),\Gam
}
\]

Let us construct a tree $\calt\subset{}^{<\ome}\ome$ recursively as follows.
For $a\in\calt$, $Seq(a)$ is a label attached with the node $a$,
which is a sequent at $a$.
First put the empty sequent at the root $\emptyset$.
{\bf Leaf condition} on the tree runs:
If $Seq(a)$ is an axiom in $\bfG(\mathcal{Q})$, then $a$ is a leaf in $\calt$. 
The construction is divided to three.
Suppose that the tree $\calt$ has been constructed up to a node $a\in{}^{<\ome}\ome$. 
\\
{\bf Case 0}. $lh(a)=3i$ for an $i\geq 0$:
Apply the inference $(W)_{i}$ backwards.
\\
{\bf Case 1}. $lh(a)=3i+1$: Apply one of inferences $(\lor),(\land),(\exi), (\fal\ome),(\exi^{2})$ if it is possible.
Otherwise repeat.
\\
{\bf Case 2}. $lh(a)=3\la n,i\ra+2$ for $n,i<\ome$:
Apply the inference $(prg)_{i}$ backwards if it is possible.
Otherwise repeat.
\\

If the tree $\calt$ is not well-founded, then let $\mathcal{P}$ be an infinite path through $\calt$.
Let $(M)_{i}\subset\Natural$ be a set such that for any $n\in\Natural$,
$(X_{i}(n))\in\mathcal{P} \Rarw n\not\in(M)_{i}$ and $(\bar{X}_{i}(n))\in\mathcal{P} \Rarw n\in(M)_{i}$.
Then for any $n$ for which one of $X_{i}(n),\bar{X}_{i}(n)$ is in $\mathcal{P}$,
we obtain $n\in(\mathcal{Q})_{i}\Lrarw n\in(M)_{i}$.
For other $n$, $n\in(M)_{i}$ is arbitrarily determined for $i\neq 0$:
set $(M)_{0}:=(\mathcal{Q})_{0}$.

$M$ is shown to be a countable coded $\ome$-model of ${\rm WOP}({\sf g})$ as follows.
The search procedure is fair, i.e., each formula is eventually analyzed on every path
as in \cite{AfshariRathjen,RathjenHFriedman,RathjenWeiermannRM}.
We see from the fairness that
$M\not\models A$ by induction on the number of occurrences of logical connectives in 
formulas $A$ on the path $\mathcal{P}$.

\section{Cut elimination}

In what follows assume that $\calt$ is well founded.
Since we are working in ${\rm ACA}_{0}$, 
we know that the Kleene-Brouwer ordering $<_{KB}$ on $\calt$ is
a well-ordering, cf.\,\cite{Simpson}.
Let $\Lam=otp(<_{KB})$ denote the order type of the well-ordering $<_{KB}$.
We have ${\rm WO}({\sf g}^{\prime}(\Lam))$ 
by ${\rm WOP}({\sf g}^{\prime})$ and ${\rm WO}(\Lam)$.

For $b<\Lam$  let us write $\vdash^{b}\Gam$ 
when there exists a derivation of $\Gam$ in $\bfG(\mathcal{Q})+(prg)+(W)$ whose depth is bounded by $b$.
On the other side for $\alp<{\sf g}^{\prime}(\Lam)$, $\vdash^{\alp}_{0}\Gam$
designates that there exists a derivation of $\Gam$ in $\bfG(\mathcal{Q})+(prg)$ of depth $\alp$.
In the derivation no inference $(W)_{i}$ occurs.
Specifically
for a function $\pi$ on ${}^{<\ome}\ome$,
$\pi\vdash^{\alp}_{0}\Gam$ designates that there exists a derivation of the sequent $\Gam$
in $\bfG(\mathcal{Q})+(prg)$ with the repetition rule $(Rep)$ in depth $\alp$, and this fact is witnessed by the function $\pi$.
The last means the following.
For each $a\in{}^{<\ome}\ome$, either $\pi(a)=*$ designating that $a$ is not in the naked tree for the derivation,
or $\pi(a)=(Seq(a), Rule(a),Mfml(a),Sfml(a),ord(a))$, where $Seq(a)$ denotes the sequent 
at the node $a$, $Rule(a)$ the inference rule whose lower sequent is $Seq(a)$,
$Mfml(a)$ is the main (principal) formula of $Rule(a)$, $Sfml(a)$ the minor (auxiliary or side) formulas
of $Rule(a)$ and $ord(a)<\Lam$.

The following Theorem \ref{th:Takeutielementaryrecursive} is due to 
G. Takeuti\cite{TRemark,PT2}\footnote{Actually Takeuti proved a similar result when we have in hand
a finite proof figure of transfinite induction in {\sf PA}.
Under the assumption we can take an order preserving map $f$ elementarily recursive in
the ordering, cf.\,\cite{attic}.}.

\bth\label{th:Takeutielementaryrecursive}
{\rm The following is provable in ${\rm RCA}_{0}+{\rm WO}(\alp)$:}
\\
Suppose that $\prec$ is a linear ordering with the least element $0$, and $<$ denotes
the well-ordering up to $\ome^{\alp}$.
$(prg)_{\prec}$
denotes the sequent calculus with inference rules $(prg)_{\prec}$ and the repetition rule $(Rep)$.
\[
\infer[(prg)_{\prec}]{\Gam, E(m)}
{
\{\Gam,E(n): n\prec m \mbox{ is true}\}
}
\]
Suppose $\pi\vdash_{0}^{\alp}\fal x\, E(x)$.
Then there exists an embedding $f$ such that $n\prec m \Rarw f(n)<f(m)$, $f(m)<\ome^{\alp}$
for any $n,m$,
and $f$ is $\alp$-recursive in the function $\pi$ and the relations $\prec, <$.
\end{theorem}
\bprf
Let us write $\Gam:\alp$ for $\vdash^{\alp}_{0}\Gam$, and $<_{\ome}$ for the usual $\ome$-ordering
 in the proof.
First search the $\ome$-rule $(\fal\ome)$ nearest to the root in the derivation $\pi$:
\[
\infer*{\fal x\, E(x):\alp}
{
\infer[(\fal\ome)]{\fal x\, E(x):\alp^{\prime}}
 {
  \{
   \infer*[\pi_{m}]{E(m):\alp_{m}}{}
   \}_{m\in\Natural}
 }
}
\]
where $\alp_{m}<\alp^{\prime}\leq\alp$ and there are some (possibly none) $(Rep)$'s below the 
$(\fal\ome)$.
Such an $(\fal\ome)$ exists by ${\rm WO}(\alp)$.
By induction on $m$, we define a derivation $\rho_{m}$ of $\Gam_{m}:\bet_{m}$ for a finite set
$\Gam_{m}\subset\{E(n):n\in\Natural\}$ such that $E(m)\in\Gam_{m}$ and
$\fal n[E(n)\in\Gam_{m}\Rarw m\preceq n]$ as follows.
If $\fal n<_{\ome}m(n\prec m)$, then $\rho_{m}=\pi_{m}$ and $\bet_{m}=\alp_{m}$.
Otherwise let
\beqn\label{eq:Takeutielementaryrecursive}
n_{0}\prec\cdots\prec n_{j-1}\prec n_{j}(=m)\prec n_{j+1}\prec\cdots n_{m}
\eeqn
with $\{n_{i}:i\leq m\}=\{0,\ldots,m\}$ and $j<_{\ome}m$.

Search the nearest inference $(prg)_{\prec}$ in $\rho_{n_{j+1}}$:
\[
\infer*[\rho_{n_{j+1}}]{\Gam_{n_{j+1}}:\bet_{n_{j+1}}}
{
 \infer[(prg)_{\prec}]{\Gam_{n_{j+1}}:\bet_{n_{j+1}}^{\prime}}
 {
  \{
   \infer*{\Gam_{n_{j+1}}, E(n):\bet}{}
   \}_{n\prec n^{\prime}}
  }
}
\]
where $\bet<\bet_{n_{j+1}}^{\prime}\leq\bet_{n_{j+1}}$, 
$E(n^{\prime})\in\Gam_{n_{j+1}}$ is the main formula of the inference $(prg)_{\prec}$.
We have $m\prec n_{j+1}\preceq n^{\prime}$.
Define $\rho_{m}$ be the following
\[
   \infer*{\Gam_{n_{j+1}}, E(m):\bet}{}
\]
with $\bet_{m}=\bet<\bet_{n_{j+1}}$.

Define a function $f(m)$ by induction on $m$ as follows.
$f(0)=\ome^{\bet_{0}}=\ome^{\alp_{0}}$ for the least element $0$ with respect to $\prec$.
For $m\neq 0$, $f(m)=f(n_{j-1})+\ome^{\bet_{m}}$ with the largest element $n_{j-1}<_{\ome}m$ 
with respect to $\prec$ in (\ref{eq:Takeutielementaryrecursive}).
Let us show that $f$ is a desired embedding.
In (\ref{eq:Takeutielementaryrecursive}), it suffices to show by induction on $m$ that
\beqn\label{eq:Takeutielementaryrecursive.1}
\fal i<_{\ome}m[f(n_{i+1})=f(n_{i})+\ome^{\bet_{n_{i+1}}}]
\eeqn
First by the definition of $f$ we have $f(m)=f(n_{j-1})+\ome^{\bet_{m}}$ with $m=n_{j}$.
On the other hand we have 
$f(m)+\ome^{\bet_{n_{j+1}}}=f(n_{j-1})+\ome^{\bet_{m}}+\ome^{\bet_{n_{j+1}}}=f(n_{j-1})+\ome^{\bet_{n_{j+1}}}=f(n_{j+1})$ by $\bet_{m}<\bet_{n_{j+1}}$ and IH.
This shows (\ref{eq:Takeutielementaryrecursive.1}), and our proof is completed.
\eprf
\\

Let us call a sequent $\Del$ an \textit{E-sequent} if
$\Del\subset\{\fal x\, E_{i}(x),E_{i}(n):i,n<\ome\}$.
An \textit{$E$-free} formula is a formula in which no $E_{i}$ occurs.

\blem\label{th:WOPelim}
For an $E$-sequent $\Del$ and an $E$-free sequent $\Gam$,
if $\vdash^{b}\Del,\Gam$ for $b<\Lam$, then
$\vdash^{{\sf g}^{\prime}(b)}_{0}\Del,\Gam$.
\elem
{\bf Proof} by induction on $b<\Lam$.
\\
{\bf Case 1}. $\Del,\Gam$ is an axiom:
There is nothing to prove.
\\
{\bf Case 2}. $\Del,\Gam$ is a lower sequent of an inference such that its 
principal formula is in $\Del\cup\Gam$:
\[
\infer{\vdash^{b}\Del,\Gam}
{
\cdots
&
\vdash^{c_{n}}\Del_{n},\Gam_{n}
&
\cdots
}
\]
By IH we have 
$\vdash^{{\sf g}^{\prime}(c_{n})}_{0}\Del_{n},\Gam_{n}$.
From ${\sf g}^{\prime}(c_{n})<{\sf g}^{\prime}(b)$ 
we obtain
\[
\infer{\vdash_{0}^{{\sf g}^{\prime}(b)}\Del,\Gam}
{
\cdots
&
\vdash_{0}^{{\sf g}^{\prime}(c_{n})}\Del_{n},\Gam_{n}
&
\cdots
}
\]
When there is no upper sequents, i.e., when $(E_{i}(m))\in\Del$ with the minimal $m$ with respect to $<_{i}$,
we have $\vdash^{0}_{0}\Del,\Gam$.
\\
{\bf Case 3}. $\Del,\Gam$ is a lower sequent of an inference $(W)_{i}$.
\[
\infer[(W)_{i}]{\vdash^{b}\Del,\Gam}
{
\vdash^{c^{\prime}}\Del,\Gam,{\rm LO}(<_{i})
&
\vdash^{c}\Del,\Gam, \fal x\,E_{i}(x)
&
\vdash^{d}\Del,\Gam,\exi Y\lnot{\rm TI}(<_{{\sf g}_{i}},Y)
}
\]
where 
$c^{\prime},c,d<b$.

If ${\rm LO}(<_{i})$ is false, i.e., $<_{(\mathcal{Q})_{i}}$ is not a linear ordering, then we see that
$\vdash^{c^{\prime}}\Del,\Gam$ with $c^{\prime}<b$.
IH yields the assertion.

In what follows assume that $<_{(\mathcal{Q})_{i}}$ is a linear ordering.
By IH we have for the $E$-sequent $\Del\cup\{\fal x\,E_{i}(x)\}$
\[
\vdash^{{\sf g}^{\prime}(c)}_{0}\Del,\fal x\,E_{i}(x),\Gam
\]

If $\vdash^{{\sf g}^{\prime}(c)}_{0}\Del,\Gam$, 
then we obtain the assertion.
Assume that this is not the case.
Then we claim that
\beqn\label{eq:Wiclm}
\vdash^{{\sf g}^{\prime}(c)}_{0}\fal x\, E_{i}(x)
\eeqn
This is seen by induction on ${\sf g}^{\prime}(c)<{\sf g}^{\prime}(\Lam)$ as follows.
If $\Del,\fal x\,E_{i}(x),\Gam$ is an axiom, then so is $\Gam$, i.e., either a true literal is in $\Gam$ or
$\{L,\bar{L}\}\subset\Gam$ for a literal $L$.
Then $\vdash^{{\sf g}^{\prime}(c)}_{0}\Del,\Gam$.
Next assume that $\Del,\fal x\,E_{i}(x),\Gam$ is derived by an inference whose principal formula is in 
$\Del\cup\Gam$.
\[
\infer{\vdash^{{\sf g}^{\prime}(c)}_{0}\Del,\fal x\,E_{i}(x),\Gam}
{
\{\vdash^{b_{n}}_{0}\Del_{n},\fal x\,E_{i}(x),\Gam_{n}\}_{n}
}
\]
We can assume that there exists an $n$ for which
 $\vdash^{b_{n}}_{0}\Del_{n},\Gam_{n}$ does not hold.
By IH we obtain $\vdash^{b_{n}}_{0}\fal x\,E_{i}(x)$.
Finally let
\[
\infer[(\fal\ome)]{\vdash^{{\sf g}^{\prime}(c)}_{0}\Del,\fal x\,E_{i}(x),\Gam}
{
\{\vdash^{b_{n}}_{0}\Del,E_{i}(n),\Gam\}_{n}
}
\]
We can assume that $\vdash^{b_{n}}_{0}\Del,\Gam$ does not hold for any $n$.
Then we show that $\vdash^{b_{n}}_{0}E_{i}(n)$ holds for any $n$ by induction on $b_{n}$.
Consider the case
\[
\infer[(prg)_{i}]{\vdash^{a}_{0}\Del,E_{i}(m),\Gam}
{
\{\vdash^{a_{n}}_{0}\Del,E_{i}(n),\Gam:n<_{i}m\}
}
\]
By IH we see that $\vdash^{a_{n}}_{0}E_{i}(n)$ for any $n<_{i}m$.
Thus (\ref{eq:Wiclm}) is shown.

Let $\bet_{0}={\sf g}^{\prime}(c)$.
By Theorem \ref{th:Takeutielementaryrecursive} 
there is an embedding $f$ such that
$n<_{i} m \Rarw f(n)<f(m)$, $f(m)<\ome^{\bet_{0}}$
for any $n,m$,
and $f$ is $\bet_{0}$-recursive in the computable function $\pi$ for the derivation witnessing the fact (\ref{eq:Wiclm}) 
and the relations $<_{i}, <$.

By Proposition \ref{prp:dilate} let $F$ be an order preserving map from ${\sf g}(<_{i})$ to $<$:
$n<_{{\sf g}_{i}}m\Rarw F(n)<F(m)$,
$F(m)<{\sf g}(\ome^{\bet_{0}})$, and $F$ is computable from $f$.

The following shows that $\vdash^{G(m)+3}_{0} \lnot Prg[<_{{\sf g}_{i}},Z], E(m)$ with a fresh variable $Z$
and $G(m)=\ome+1+4F(m)$
by induction on $F(m)$: 
{\small
\[
\infer[(\exi)]{\vdash^{G(m)+3}_{0}\lnot Prg[<_{{\sf g}_{i}},Z] , Z(m)}
{
\infer[(\land)]{\vdash^{G(m)+2}_{0}\lnot Prg[<_{{\sf g}_{i}},Z], \fal y<_{{\sf g}_{i}}m\, Z(y)\land \bar{Z}(m), Z(m)}
{
 \infer[(\fal\ome)]{\vdash^{G(m)+1}_{0} \lnot Prg[<_{{\sf g}_{i}},Z], \fal y<_{{\sf g}_{i}}m\, Z(y)}
 {
  \infer[(\lor)]{\{\vdash^{G(m)}_{0} \lnot Prg[<_{{\sf g}_{i}},Z], (n\not<_{{\sf g}_{i}}m)\lor Z(n) :n\in\ome\}}
  {
   \{\vdash^{G(n)+3}_{0}\lnot Prg[<_{{\sf g}_{i}},Z],  Z(n): n<_{{\sf g}_{i}}m\}
   &
   \{\vdash^{\ome}_{0} n\not<_{{\sf g}_{i}}m : n\not<_{{\sf g}_{i}}m\}
   }
  }
 &
 \hskip -10mm
 \vdash^{0}_{0}\bar{Z}(m), Z(m)
 }
}
\]
}
where $n\not<_{{\sf g}_{i}}m$ denotes the formula $\lnot({\sf g}(X_{i}))(\la n,m\ra)$,
which is a $\Del_{1}$-formula in $X_{i}$.
Thus for $G(m)+3<{\sf g}(\ome^{\bet_{0}})$ we obtain
\[
\infer[(\fal^{2})]{\vdash^{{\sf g}(\ome^{\bet_{0}})+3}_{0}  \fal Y{\rm TI}(<_{{\sf g}_{i}},Y)}
{
\infer[(\lor)]{\vdash^{{\sf g}(\ome^{\bet_{0}})+2}_{0} {\rm TI}(<_{{\sf g}_{i}},Z)}
{
 \infer[(\fal\ome)]{\vdash^{{\sf g}(\ome^{\bet_{0}})}_{0} \lnot Prg[<_{{\sf g}_{i}},Z], \fal x\, Z(x)}
 {
 \{\vdash^{G(m)+3}_{0}\lnot Prg[<_{{\sf g}_{i}},Z], Z(m):m<\ome\}
 }
}
}
\]
On the other hand we have by IH
$\vdash^{{\sf g}^{\prime}(d)}_{0}\Del,\Gam,\exi Y\lnot{\rm TI}(<_{{\sf g}_{i}},Y)$
for the $E$-free sequent $\Gam\cup\{\exi Y\lnot{\rm TI}(<_{{\sf g}_{i}},Y)\}$.
By cut-elimination we obtain
$\vdash_{0}^{\bet_{1}}\Del,\Gam$ for 
$\bet_{1}=\ome_{k}({\sf g}(\ome^{\bet_{0}})\#3\#{\sf g}^{\prime}(d))$
for a $k<\ome$ depending only on the formula $ \fal Y{\rm TI}(<_{{\sf g}_{i}},Y)$.
Now
$\bet_{1}=\ome_{k}({\sf g}(\ome^{{\sf g}^{\prime}(c)})\#3\#{\sf g}^{\prime}(d))
<{\sf g}^{\prime}(b)$
since $c,d<b$ and ${\sf g}^{\prime}(b)$ is closed under $+,\ome$ and ${\sf g}$ by (\ref{df:omegag.4}).
\eprf
\\

Let us finish our proof of the harder direction in Theorem \ref{th:derivativemodel}.
By our assumption we have $\vdash^{b}\emptyset$ for a $b<\Lam$ and the empty sequent 
$\emptyset$.
Lemma \ref{th:WOPelim} yields $\vdash^{{\sf g}^{\prime}(b)}_{0}\emptyset$.
We see that this is not the case by induction on ${\sf g}^{\prime}(b)<{\sf g}^{\prime}(\Lam)$.
Therefore the tree $\calt$ is not well founded.

Finally let us spend a few words on a formalization of the above proof in ${\rm ACA}_{0}$.
In the proof one can agree that each infinite derivation is a computable function $\pi$
on the set of finite sequences $a$ of natural numbers.
$\pi(a)$ is a bunch of data as described before Theorem \ref{th:Takeutielementaryrecursive}.
$\vdash^{\alp}\Gam$ denotes the fact that there exists a computable function $\pi$ such that
$Seq(\emptyset)=\Gam$ and $ord(\emptyset)=\alp$ for the empty sequence $\emptyset$, i.e.,
the root of the derivation tree.
$\vdash^{\alp}\Gam$ is arithmetically definable, defined by a $\Sig^{0}_{3}$-formula.
The above proof of Lemma \ref{th:WOPelim} is formalizable in ${\rm ACA}_{0}$
with the assumption ${\rm WO}({\sf g}^{\prime}(\Lam))$.

\brem
{\rm
We can show one of equivalences due to Girard\cite{Girard} in the spirit of Rathjen\cite{AfshariRathjen,RathjenHFriedman,RathjenWeiermannRM}:
${\rm ACA}_{0}$ is equivalent to ${\rm WFP}(\lam X. 2^{X})$ over ${\rm RCA}_{0}$,
where 
${\rm WFP}({\sf g}) :\Lrarw \fal X\left( {\rm WF}(X)\to {\rm WF}({\sf g}(X)) \right)$  with
${\rm WF}(X)  :\Lrarw  \fal Y\, {\rm TI}(<_{X},Y)$.

The direction ${\rm ACA}_{0}\to{\rm WFP}(\lam X. 2^{X})$ is well known.
The reverse direction is seen as follows.
Consider the proof search of the contradiction in a sequent calculus $\bfG(\mathcal{Q})+(Jcut)+(J)$.
Pick a fresh unary predicate symbol $J$. Let $\exi x\, B(x,y)$ be a fixed $\Sig_{1}$-formula.
$\bfG(\mathcal{Q})+(Jcut)+(J)$ is obtained from the sequent calculus $\bfG(\mathcal{Q})$ by adding the following three
inference rules $(Jcut),(J),(\bar{J})$.
\[
\infer[(Jcut)]{\Gam}
{
\Gam,J(n)
&
\bar{J}(n),\Gam
}
\,
\infer[(J)]{\Gam,J(n)}
{
\Gam,J(n),\exi x\, B(x,n)
}
\,
\infer[(\bar{J})]{\Gam,\bar{J}(n)}
{
\Gam,\bar{J}(n),\fal x\, \lnot B(x,n)
}
\]
$J(n)$ is intended to denote $\exi x\, B(x,n)$.
If the tree in the proof search is not well founded, then an infinite path through the tree yields a set $\mathcal{J}$ such that
$\fal n[n\in \mathcal{J}\lrarw \exi x\, B(x,n)]$. Thus ${\rm ACA}_{0}$ follows.
Suppose contrarily that the tree is well founded, and let $\Lam$ be the depth of the well founded tree.
Then a cut elimination yields a cut-free derivation of the empty sequent in $\bfG(\mathcal{Q})$ in depth
$2_{c}(\Lam)$ for a constant $c$ depending only on the $\Del_{0}$-formula $B$.
From ${\rm WFP}(\lam X. 2^{X})$ we see that the cut-free derivation is well founded, and this is not the case.
}
\erem

\section{Common fixed points}
Let $\alp$ be the order type of a computable well ordering on $\Natural$.
$\vphi[{\sf g}]_{\alp}(\bet)$ denotes the $\alp$-th Veblen function starting with 
$\vphi[{\sf g}]_{0}\bet={\sf g}(\bet)$.

We assume that $\vphi[{\sf g}]_{\alp}(X)$ is a term structure over constants $\{\vphi[{\sf g}]_{\alp}( c ) :c\in X\cup\{0\}\}$ and
unary function symbols $\vphi[{\sf g}]_{\bet}\,(\bet<\alp)$ and the addition $+$.
Also a function symbol for the exponential $\ome^{\bet}$ is included when
$\vphi[{\sf g}]_{0}={\sf g}$ is not the exponential.

In what follows we assume that
each term structure $\vphi[{\sf g}]_{\bet}(X)\, (\bet<\alp)$ is extendible,
and $\vphi[{\sf g}]_{\alp}$ is exponential.
Moreover we assume that the followings are provable in ${\rm RCA}_{0}$, cf.\,(\ref{df:omegag.4}).

\beqnarr
\bet_{1},\ldots,\bet_{n}<_{\vphi[{\sf g}]_{\alp}(X)}\vphi[{\sf g}]_{\alp}(c) & \to &
 f(\bet_{1},\ldots,\bet_{n})<_{\vphi[{\sf g}]_{\alp}(X)}\vphi[{\sf g}]_{\alp}(c) 
 \nonumber
 \\
&& (f\in\{\vphi[{\sf g}]_{\bet}:\bet<\alp\}\cup\{+\})
 \nonumber
 \\
\ome^{\vphi[{\sf g}]_{\alp}(c)} & = & \vphi[{\sf g}]_{\bet}(\vphi[{\sf g}]_{\alp}(c))=\vphi[{\sf g}]_{\alp}(c)\, (\bet<\alp)
\nonumber
\\
\vphi[{\sf g}]_{\alp}(0) & = & \sup\{(\vphi[{\sf g}]_{\bet})^{n}(0):\bet<\alp, n\in\ome\}
\label{eq:Veblenfund}
\\
\vphi[{\sf g}]_{\alp}(c+1) & = & \sup\{(\vphi[{\sf g}]_{\bet})^{n}(\vphi[{\sf g}]_{\alp}(c)+1):\bet<\alp, n\in\ome\}
\nonumber
\eeqnarr

We see the following as in Proposition \ref{prp:dilate}.

\bprp\label{prp:dilatevphi}
{\rm Suppose $\vphi[{\sf g}]_{\bet}(X)$ is an extendible term structure.
Then the following is provable in ${\rm RCA}_{0}$:}
Let both $X$ and $Y$ be linear orderings.

Let $f:\{0\}\cup X\to \{0\}\cup Y$ be an order preserving map, $n<_{X}m\Rarw f(n)<_{Y}f(m)\, (n,m\in\{0\}\cup X)$.
Then there is an order preserving map $F:\vphi[{\sf g}]_{\bet}(X)\to\vphi[{\sf g}]_{\bet}(Y)$, 
$n<_{\vphi[{\sf g}]_{\bet}(X)}m\Rarw F(n)<_{\vphi[{\sf g}]_{\bet}(Y)}F(m)$.
\eprp

\bth\label{th:commonderivative}
Let each term structure $\vphi[{\sf g}]_{\bet}(X)\, (\bet<\alp)$ be extendible,
and $\vphi[{\sf g}]_{\alp}$ is exponential for which (\ref{eq:Veblenfund}) holds.
Then the following two are mutually equivalent over ${\rm ACA}_{0}+{\rm LO}(\alp)$
for $\alp>0$.
\benu
\item
${\rm WOP}(\vphi[{\sf g}]_{\alp})$.
\item
$\left(\fal\bet<\alp{\rm WOP}(\vphi[{\sf g}]_{\bet})\right)^{+}$.
\eenu
\end{theorem}

The easier direction states that ${\rm WOP}(\vphi[{\sf g}]_{\alp})$ follows from
$\left(\fal\bet<\alp{\rm WOP}(\vphi[{\sf g}]_{\bet})\right)^{+}$, and
it follows from the fact (\ref{eq:Veblenfund}).

The harder direction is seen as in Theorem \ref{th:derivativemodel} by slight modifications.
Suppose ${\rm WOP}(\vphi[{\sf g}]_{\alp})$ and ${\rm LO}(\alp)$ for $\alp>0$.
Replace the inference rule $(W)_{i}$ by
\[
\infer[(W)_{i,\bet}]{\Gam}
{
\Gam,{\rm LO}(<_{i})
&
\Gam,\fal x\,E_{i}(x)
&
\exi Y\lnot{\rm TI}(<_{i,\bet},Y),\Gam
}
\]
where $\bet<\alp$ and $n<_{i,\bet}m:\Lrarw n<_{\vphi[{\sf g}]_{\bet}(<_{i})}m$.

Construct fairly a tree in the sequent calculus $\bfG(\mathcal{Q})+(prg)+\{(W)_{\bet}\}_{\bet<\alp}$
ending with the empty sequent.
When the tree is not well founded, an infinite path through the tree yields a countable coded $\ome$-model
of $\fal\bet<\alp{\rm WOP}(\vphi[{\sf g}]_{\bet})$.

Suppose that the search tree $\calt$ is well founded with its order type $\Lam$ 
in the Kleene-Brouwer ordering.
We obtain ${\rm WO}(\vphi[{\sf g}]_{\alp}(\Lam))$ by ${\rm WOP}(\vphi[{\sf g}]_{\alp})$.
As in Lemma \ref{th:WOPelim} we see the following Lemma \ref{th:WOPelimvphi}
from Proposition \ref{prp:dilatevphi} and (\ref{eq:Veblenfund}):
$c,d<b \spand \bet<\alp \Rarw 
\ome_{k}(\vphi[{\sf g}]_{\bet}(\ome^{\vphi[{\sf g}]_{\alp}(c)})\# 3\#\vphi[{\sf g}]_{\alp}(d))<\vphi[{\sf g}]_{\alp}(b)$

\blem\label{th:WOPelimvphi}
For an $E$-sequent $\Del$ and an $E$-free sequent $\Gam$,
if $\vdash^{b}\Del,\Gam$ for $b<\Lam$, then
$\vdash^{\vphi[{\sf g}]_{\alp}(b)}_{0}\Del,\Gam$.
\elem

The harder direction in Theorem \ref{th:commonderivative} is concluded as follows.
By our assumption we have $\vdash^{b}\emptyset$ for a $b<\Lam$ and the empty sequent 
$\emptyset$.
Lemma \ref{th:WOPelimvphi} yields $\vdash^{\vphi[{\sf g}]_{\alp}(b)}_{0}\emptyset$.
We see that this is not the case by induction on $\vphi[{\sf g}]_{\alp}(b)<\vphi[{\sf g}]_{\alp}(\Lam)$.
Therefore the tree $\calt$ is not well founded.



\end{document}